\documentclass[11pt,leqno]{article}
\usepackage{amssymb,amsmath,amsthm}
\usepackage[all]{xy}
\CompileMatrices
\numberwithin{equation}{section}

\newcommand{\rmpt}{{\rm pt}}

\newcommand{\id}{{\rm id}}

\def\etens{\mathbin{\boxtimes}}
\newcount\temp
\temp=\catcode`\@
\catcode`\@=11
\def\@bletens{\mathbin{\etens^{L}}}
\def\@letens_#1{\mathbin{\etens_{\raise1.5ex\hbox to-.1em{}#1}^{L}}}
\def\letens{\@ifnextchar _{\@letens}{\@bletens}}
\catcode`\@=\temp

\def\rop{{\rm op}}

\def\phi{{\varphi}}
\def\epsilon{\varepsilon}

\newcommand{\mds[1]}{{\mathfrak{Mod}}({#1})}

\newcommand{\mdso[1]}{{\mathfrak{Mod}}_{0}({#1})}

\newcommand{\md[1]}{{\rm Mod}({#1})}

\def\dim{{\rm dim}}

\newcommand{\stx}{\operatorname{{\mathfrak{X}}}}
\newcommand{\sty}{\operatorname{{\mathfrak{Y}}}}

\newcommand{\cor}{\operatorname{\mathbf{k}}}

\newcommand{\wcor}{\widehat\cor}

\newcommand{\ba}{\begin{array}}
\newcommand{\ea}{\end{array}}

\def\sha{\mathcal{A}}

\def\shd{\mathcal{D}}
\def\she{\mathcal{E}}
\def\shf{\mathcal{F}}

\def\shi{\mathcal{I}}

\def\shl{\mathcal{L}}

\def\sho{\mathcal{O}}

\def\shr{\mathcal{R}}
\def\shs{\mathcal{S}}
\def\sht{\mathcal{T}}
\def\shu{\mathcal{U}}
\def\shv{\mathcal{V}}
\def\shw{\mathcal{W}}

\def\shde{\she^{\sqrt{v}}}
\def\shhe{\widehat{\she}}

\def\shie{\she^{\infty}}
\def\shiw{\shw^{\infty}}

\newcommand{\shed[1]}{\shde_{#1}}
\newcommand{\shhed[1]}{\shhde_{#1}}

\def\shdw{\shw^{\sqrt{v}}}
\def\shhw{\widehat{\shw}}

\newcommand{\shwd[1]}{\shdw_{#1}}
\newcommand{\shhwd[1]}{\shhdw_{#1}}

\newcommand{\shetd[1]}{\she^{\sqrt{v}}_{#1\times\C,\hat{t}}}
\newcommand{\shet[1]}{\she_{#1\times\C,\hat{t}}}

\newcommand{\hcor}{\widehat\cor}

\newcommand{\Ad}{\mathrm{Ad}}

\newcommand{\stkHom}[1][]{\mathfrak{Hom}_{\raise1.5ex\hbox to.1em{}#1}}

\newcommand{\lienR}{\mathsf{R}}
\newcommand{\lienS}{\mathsf{S}}
\newcommand{\lienT}{\mathsf{T}}

\newcommand{\lien}{\mathsf}

\newcommand{\C}{\mathbb{C}}
\newcommand{\N}{\mathbb{N}}

\newcommand{\Z}{\mathbb{Z}}

\renewcommand{\to}[1][]{\xrightarrow[#1]{}}

\newcommand{\isoto}[1][]{\xrightarrow[#1]{\sim}}

\newcommand{\tens}{\otimes}
\newcommand{\eim}[1]{{#1}_{!}}
\newcommand{\oim}[1]{{#1}_*}
\newcommand{\opb}[1]{#1^{-1}}
\renewcommand{\hom}[1][]{{\mathcal{H}om}_{\raise1.5ex\hbox to.1em{}#1}}

\newcommand{\rhom}[1][]{{R\mathcal{H}om}_{\raise1.5ex\hbox to.1em{}#1}}

\newcommand{\sect}{\Gamma}

\newcommand{\Hom}[1][]{\mathrm{Hom}_{\raise1.5ex\hbox to.1em{}#1}}
\newcommand{\RHom}[1][]{\mathrm{RHom}_{\raise1.5ex\hbox to.1em{}#1}}

\theoremstyle{plain}

\newtheorem{theorem}{Theorem}[section]
\newtheorem{proposition}[theorem]{Proposition}
\newtheorem{lemma}[theorem]{Lemma}

\theoremstyle{definition}

\newtheorem{definition}[theorem]{Definition}
\newtheorem{notation}[theorem]{Notation}

\newtheorem{remark}[theorem]{Remark}
\newtheorem{remarks}[theorem]{Remarks}

\newcommand{\eq}{\begin{eqnarray}}
\newcommand{\eneq}{\end{eqnarray}}
\newcommand{\eqn}{\begin{eqnarray*}}
\newcommand{\eneqn}{\end{eqnarray*}}

\newcommand{\lp}{{\rm(}}
\newcommand{\rp}{{\rm)}}

\newenvironment{nnum}{

  \begin{enumerate}
  \itemsep=0pt
  
  }
  {\end{enumerate}}

\newenvironment{anum}{
  \begin{enumerate}
  \itemsep=0pt
  
  }
  {\end{enumerate}}
\newcommand{\bnum}{\begin{nnum}}
\newcommand{\enum}{\end{nnum}}
\newcommand{\banum}{\begin{anum}}
\newcommand{\eanum}{\end{anum}}

\def\dddt{{\raise-.3em\hbox{$\big\cdot$}}}

\author{Pietro Polesello and Pierre Schapira}

\title{Stacks of quantization-deformation modules on
  complex symplectic manifolds}

\begin{document}

\maketitle

\begin{abstract}
On a complex symplectic manifold $\stx$, we
construct the stack of quantization-deformation modules, that is, 
(twisted) modules of microdifferential operators with an extra 
central parameter $\tau$, a substitute to the lack of homogeneity.
We also quantize involutive submanifolds of contact manifolds.
\end{abstract}

\tableofcontents

\section{Introduction}
Masaki Kashiwara \cite{K2} has constructed the stack of modules of 
microdifferential 
operators on a complex contact manifold $\sty$. Let us explain
briefly what this means. A local model for $\sty$ 
is an open subset $V$ of $P^*X$, the projective cotangent bundle of a
complex manifold $X$. The manifold $P^*X$ is endowed with 
the sheaf of rings $\she_X$ of microdifferential operators of 
Sato-Kawai-Kashiwara \cite{S-K-K}. Twisting this ring by half-densities on 
$X$, we find another ring $\shed[X]$ with an extra property: namely, it 
is endowed with an anti-involution $*$.
If $\psi:P^*X\supset V_X\isoto V_Y\subset P^*Y$ is a contact transformation,  
one can locally ``quantize'' it as a ring isomorphism 
$\Psi: \oim{\psi}(\shed[X]\vert_{V_X})\isoto (\shed[Y]\vert_{V_Y})$ 
commuting with the anti-involution $*$. This quantization is not
unique and this fact makes it impossible to glue together the sheaves
$\shed[X]\vert_V$'s and to get a globally defined sheaf of rings on
$\sty$. However, and this is the content of Kashiwara's paper, it is
possible to glue together the categories of abelian sheaves
$\md[{\shed[X]\vert_V}]$  and to get a canonically defined stack 
(roughly speaking, a sheaf of categories) on $\sty$. 
The proof has two aspects, one analytical, 
which consists essentially in noticing that the automorphisms of the
ring $\shed[X]$ preserving $*$ are in bijection with a subgroup of its
invertible elements, the  
other one purely algebraic, dealing with the machinery of stacks.

In this paper, we first recall Kashiwara's proof and extend it to the case
of regular involutive submanifolds of $\sty$. Then, and this is our
main result, we adapt it to the case of symplectic manifolds. As we
shall see, new difficulties appear.

The local model is now an open subset $U$ of the cotangent
bundle $T^*X$ of a complex manifold $X$. Let $\C$ denote the
complex line with holomorphic coordinate $t$, let $\dot{T}^*\C$ denote
the cotangent bundle to $\C$ with the zero-section removed and
with coordinates $(t,\tau)$.
Denote by $\dot{P}^*(X\times\C)$ the quotient of 
$T^*X\times\dot{T}^*\C$ by the diagonal $\C^\times$-action. Then 
$\dot{P}^*(X\times\C)$  is an open subset of the projective
cotangent bundle ${P}^*(X\times\C)$ and there is a natural 
map $\rho:\dot{P}^*(X\times\C)\to T^*X$, $(p,(t;\tau))\mapsto p \opb \tau$. 
The manifold $\dot P^*(X\times\C)$ is endowed with the sheaf of rings
$\shet[X]$ of microdifferential operators on $P^*(X\times\C)$
which commute with $D_t$ ({\em i.e.}, which do not depend on
the $t$-variable). We endow $T^*X$ with the sheaf of rings 
$\shw_X:=\oim{\rho}\shet[X]$. Roughly speaking, 
$\shw_X$ is the sheaf of microdifferential operators in the
$(x,D_x)$-variables and a central extra parameter $\tau$ of order
$1$, which kills the homogeneity of $T^*X$. 
Such algebras in the formal  case over real symplectic manifolds
are called semi-classical star-algebras, or also,
quantization-deformation algebras by many authors and 
we refer to \cite{Bo} for their study. 
Note that the link between the sheaves of rings $\she_X$ and 
$\shw_X$ is well-known from the specialists and appears explicitely 
when $\dim X=1$ in \cite{AKKT}. By reference to the WKB-method of the
physicists, these authors 
call WKB-differential operators the sections of $\shw_X$ and we shall
follow this terminology.

Denote by 
$\wcor:=\C[\tau,\opb{\tau}]\!]$ the field 
of formal Laurent series $\sum_{j\in\Z}a_j{\tau}^j$  
with $a_j=0$ for $j\gg 0$, and 
by $\cor$ the field $\shw_{\rmpt}$, a subfield of $\wcor$
(see Definition \ref{def:cor}). 

The sheaf $\shw_X$ is 
thus a sheaf of $\cor_{T^*X}$-central algebras and 
the center of $\shw_X$ is now too large in order to apply
Kashiwara's method. We overcome this difficulty 
by showing that above a symplectic transformation 
$\phi:T^*X\supset U_X\isoto U_Y\subset T^*Y$, there exists locally a contact 
transformation $\psi: \dot{P}^*(X\times\C)\supset \opb\rho (U_X)\isoto
\opb \rho (U_Y)\subset \dot{P}^*(Y\times\C)$ commuting with $\tau$ and that 
this transformation may be quantized as an isomorphism of rings 
$\Phi:\oim{\rho}\oim{\psi}(\shed[X\times\C]|_{\opb\rho(U_X)})
\isoto\oim{\rho}(\shed[Y\times\C]|_{\opb\rho (U_Y)})$, this
isomorphism $\Phi$ commuting with $D_t$, hence interchanging 
$\shwd[X]$ and $\shwd[Y]$. In general, these isomorphisms do not
allow us to glue the categories 
$\md[{\oim{\rho}(\shed[X\times\C]|_{\opb\rho(U_X)})}]$ and to obtain a
stack, due to a kind of translation operator which appears in the
fibers of $\rho$. But, fortunately, these translations act trivialy on 
$\shwd[X]$, and we can glue the categories
$\md[{\shwd[X|_{U_X}]}]$.

Note that Maxim Kontsevich \cite{Ko}  
has recently announced a similar result in the
much more general setting of (algebraic)
Poisson manifolds, based on a different method.
Also note that semi-classical star-algebras on complex symplectic 
manifolds are constructed 
under suitable hypotheses in \cite{N-T}, using Fedosov connections.
We refer to \cite{B-S} for a discussion on possible physical applications of 
such constructions.

The contents of this paper is as follows.

In Section \ref{section:stacks} we recall some facts (well
known from a few specialists) concerning the construction of
stacks. References are made to \cite{Gi}, \cite{Br}, \cite{K2}, \cite{K-S2}.

In Section \ref{section:geometry} we explain how to construct 
$\tau$-preserving contact isomorphisms associated with symplectic
isomorphisms.

In Section \ref{section:Emod} we overview the theory of
Sato-Kawai-Kashiwara \cite{S-K-K}  of 
microdifferential operators and modules over such rings 
(see also \cite{S},  \cite{K3} for a detailed exposition). 

In Section \ref{section:Cquant} we give a detailed proof of 
Kashiwara's quantization theorem of \cite{K2}.

In Section \ref{section:Invquant} we treat the case of involutive 
submanifolds of complex contact manifolds.

In Section \ref{section:Wmod} we construct the ring $\shw_X$ of
WKB-differential operators on the cotangent bundle to a complex
manifold $X$.

In Section \ref{section:QCSM} we adapt Kashiwara's proof of \cite{K2}
to construct the stack of WKB-modules on a complex symplectic manifold.

All our results extend to the formal case,  that is, to the ring 
$\shhw$ associated with the ring of formal microdifferential
operators $\shhe$.
\medskip

\noindent
{\bf Aknowledgement.}
We have benefited from extremely valuable
suggestions and detailed explanations of Masaki Kashiwara, especially
concerning Sections \ref{section:geometry}--\ref{section:Emod}. 
It is a pleasure to thank him here. We would like to thank also
Louis Boutet de Monvel for fruitful discussions. 

\section{Notations}\label{section:notation}

We will mainly follow the notations of \cite{K-S1}.
In this paper, unless otherwise specified,  all manifolds are
complex analytic. 

Let $X$ be a manifold. 
We denote by $T^*X$ the cotangent bundle, $\pi:T^*X\to X$ the
projection, $\dot{T}^*X$ the bundle $T^*X\setminus X$, where $X$ is
identified with the zero-section, $a:T^*X\to T^*X$ the antipodal map.
We shall also consider the projective cotangent bundle,
$P^*X=\dot{T}^*X/\C^\times$, ($\C^\times$ denotes the
multiplicative group of non-zero complex numbers).
We keep the notation $\pi$ for the projection $P^*X\to X$.

On a complex manifold $X$, we consider the structure sheaf  $\sho_X$, 
the sheaves  $\Omega^p_X$ of holomorphic
$p$-forms and the sheaf $\shd_X$ of linear 
holomorphic differential operators of finite order. One sets 
$\Omega_X:=\Omega^n_X$, where $n$ is the dimension of $X$.

This paper deals with stacks, which roughly speaking, means sheaves of 
categories. The classical reference is 
\cite{Gi} and a more popular one is \cite{Bry}. This notion is
also well explained in \cite{K2}.

Recall that if $\shr$ is a sheaf of unital rings on $X$, one denotes by 
$\shr^\times$ the sheaf of invertible sections, and if 
$a\in\shr^\times$ is a local section, one defines the ring automorphism 
$\Ad(a)$ of $\shr$ by setting for any local section $b\in\shr$
\eqn
&& \Ad(a)(b) = ab\opb a.
\eneqn
By an $\shr$-module (resp., $\shr^\rop$-module), we mean a sheaf of left 
(resp., right) $\shr$-modules. We denote by $\md[\shr]$ the abelian category 
of $\shr$-modules and by $\mds[\shr]$ the corresponding abelian stack 
on $X$ given by the assignement $X\supset U \mapsto \md[\shr\mid_U]$.

\section{Construction of stacks and functors}\label{section:stacks}

Let $X$ be a topological space and let $\shu=\{U_i\}_{i\in I}$ be an open
covering of $X$. Define:
\eqn
X_0=\bigsqcup_{i\in I}U_i,\quad
X_1= \bigsqcup_{i,j\in I}U_{ij}, \quad
X_2=\bigsqcup_{i,j,k\in I}U_{ijk},\mbox{ etc.}
\eneqn
where $U_{ij}=U_i\cap U_j$, $U_{ijk}=U_i\cap U_j\cap U_k$, etc.
In other words,
\eqn
X_1:={X_0}\times_X{X_0},\quad
X_n=X_0\times_X\cdots\times_X X_0 \mbox{ ($n+1$-times)}.
\eneqn

We  introduce the projections 
\eqn
&& p^n:X_n\to X,  \quad 0\leq n, \\
&& p^n_{i}:X_n\to X_0,  \quad 0\leq i\leq n,\\
&& p^n_{ij}:X_n\to X_1 ,  \quad 0\leq i<j\leq n,\\
&& p^n_{ijk}:X_n\to X_2 ,  \quad 0\leq i<j<k\leq n, \mbox{ etc.}
\eneqn
For example, $p^n_{ij}:X_n\to X_1 $ is the projection to the
$(i,j)$-factor.
In the sequel, we set $p=p^0$ and, if there is no risk of confusion, we write
$p_{i}$ instead of $p^n_{i}$, $p_{ij}$ instead of $p^n_{ij}$, etc.
Hence we have the following diagram
\eq\label{diag1}
\xymatrix{ X_{\bullet}\colon \ar[d]|-{p} 
\ar@{.>}@<-1.8ex>[r] \ar@{.>}@<-.6ex>[r] \ar@{.>}@<.6ex>[r]\ar@{.>}@<1.8ex>[r]
&X_2 
\ar@<-1.5ex>|-{p_{02}}[rr]\ar[rr]|-{p_{12}}\ar@<1.5ex>[rr]|-{p_{01}}
                                        \ar@<-1.5ex>[drrrr]|-{p^2} &&
X_1 \ar@<-1.ex>[rr]|-{p_{1}}\ar@<1.ex>[rr]|-{p_{0}} \ar@<-1.ex>[drr]|-{p^1}&&
X_0 \ar[d]|-{p} \\
X &&&&& X .}
\eneq

\begin{notation}
(i) For a sheaf $\shf$ or a stack $\mathfrak{S}$ on $X_0$, we set \
$\shf_i=\opb{p_i}\shf$ and $\mathfrak{S}_i=\opb{p_i}\mathfrak{S}$ 
for $0\leq i\leq n$, where $p_i:X_n\to X_0$, $n\geq 1$.

\noindent
(ii) For a morphism of sheaves $f$ or a functor of stacks $F$ on $X_1$, we set 
$f_{ij}=\opb{p_{ij}}f$ and $F_{ij}=\opb{p_{ij}}F$
for $0\leq i<j\leq n$, where $p_{ij}:X_n\to X_1$, $n\geq 2$.

\noindent
(iii) For a section of a sheaf $a$ or an isomorphim of functors $\alpha$ 
on $X_2$, we set $a_{ijk}=\opb{p_{ijk}}a$ and $\alpha_{ijk} = 
\opb{p_{ijk}}\alpha$ for $0\leq i<j<k\leq n$, where $p_{ijk}:X_n\to X_2$, 
$n\geq 3$. 
\end{notation}

Let $\kappa$ be a sheaf of commutative unital rings on $X$.

\begin{definition}\label{def:lien1}
A $\kappa$-lien on $X_{\bullet}\to X$ is a triplet 
$\lienR=(\shr,f,a)$,
where $\shr$ is a sheaf  of central $\opb p\kappa$-algebras on $X_0$ 
({\em i.e.,} $\opb p\kappa$ is the center of $\shr$),
$f:\shr_1\isoto\shr_0$ is an isomorphism of $\kappa_0$-algebras 
on $X_1$ and  $a$ is section in $\sect(X_2,\shr_0^\times)$, such that on $X_2$
\eq\label{eq:liens}
&& f_{01}\circ f_{12}=\Ad(a)\circ f_{02}.
\eneq
\end{definition}

On $X_3$ one has 
\eqn
({f}_{01}\circ {f}_{12})\circ {f}_{23}
       =\Ad(a_{012})\circ {f}_{02}\circ {f}_{23}
       =\Ad(a_{012}a_{023})\circ {f}_{03},
\eneqn
and 
\eqn
&&{f}_{01}\circ ({f}_{12}\circ{f}_{23}) ={f}_{01}
\circ \Ad(a_{123})\circ {f}_{13} \\
&&\hspace{3cm} = \Ad({f}_{01}(a_{123}))\circ {f}_{01}\circ {f}_{13} = 
\Ad({f}_{01}(a_{123})a_{013})\circ {f}_{03}.
\eneqn
It follows that there exists a section  $c \in \sect(X_3;\kappa^\times)$, 
such that 
\eq\label{eq:c}
&&a_{012}a_{023} = {f}_{01}(a_{123})a_{013}\cdot c \quad \mbox{ in }
\sect(X_3;\shr_0^\times).
\eneq

\begin{definition}\label{def:lien2}
One says that a $\kappa$-lien $\lienR=(\shr,f,a)$ is effective if $c=1$.
\end{definition}

\begin{notation}
Let  $\lienR=(\shr,f,a)$ be a $\kappa$-lien on $X_0$. 
We denote by $_{\opb f}(\cdot)$ the equivalence of stacks $\mds[\shr_1]\isoto
\mds[\shr_0]$ which associates to an $\shr_1$-module $\shf$ on $X_1$, the 
$\shr_0$-module $_{\opb f}\shf$, {\em i.e.} the $\opb p\kappa$-module $\shf$ 
endowed with the action induced by $\opb f\colon\shr_0\isoto\shr_1$ : 
if $l\in \shf$ and $r\in\shr_0$ are local sections, the action of $r$ on $l$ 
is given by $(r,l)\mapsto\opb f(r) l$.
 
\end{notation}

\begin{theorem}{\rm (M. Kashiwara \,\cite{K2}.)}\label{th:elienK}
To an effective $\kappa$-lien $\lienR=(\shr,f,a)$ on $X_{\bullet}\to X$ one 
associates an abelian $\kappa$-stack $\mds[\lienR]$ on $X$, an equivalence of 
$\opb p\kappa$-stacks $F_{\lienR}\colon\opb p \mds[\lienR]\isoto \mds[\shr]$ 
and an isomorphism of functors $\alpha_{\lienR}\colon _{\opb f}(\cdot) \circ 
(F_{\lienR})_1 \isoto (F_{\lienR})_0$ such that $(\alpha_{\lienR})_{01}\circ
(\alpha_{\lienR})_{12} = (\alpha_{\lienR})_{02}\circ (a \cdot)$ 
(here $a\cdot$ denotes the isomorphism of functors defined by left 
multiplication by $a$).
Moreover, the datum of $(\mds[\lienR],F_{\lienR},\alpha_{\lienR})$ is unique 
up to equivalence where the equivalence is unique up to unique isomorphism.
\end{theorem}


\begin{proof}[Sketch of the proof]
For $V$ open in $X_0$, let $\mdso[\lienR](V)$ be the category defined as 
follows: an object of $\mdso[\lienR](V)$ is a pair $(\shf,m)$ where $\shf$ is 
an $\shr$-module on $V$ and $m\colon _{\opb f}\shf_1\isoto \shf_0$ is an 
isomorphism of $\shr_0$-modules on $V_1=V\times_{X} X_1$ 
such that the following
diagram of isomorphisms of $\shr_0$-modules on $V_2=V\times_{X} X_2$
commutes
$$\xymatrix{
_{\opb f_{01}}(_{\opb f_{12}}\shf_2)\ar[r]^-{a \cdot }\ar[d]^-{m_{12}} & 
_{\opb f_{02}}\shf_2\ar[d]^-{m_{02}}\\
_{\opb f_{01}}\shf_1\ar[r]^-{m_{01}} & {\shf_0} .}
$$
Here $a\cdot$ denotes left multiplication by $a$ in $_{\opb f_{02}}\shf_2$.
Morphisms $\alpha\colon (\shf,m)\to (\shf',m')$ are morphisms of $\shr$-modules
on $V$, $\alpha\colon \shf\to \shf'$ such that the following diagram of 
morphisms of $\shr_0$-modules on $V_1$ commutes 
$$\xymatrix{
_{\opb f}\shf_1\ar[r]^-{m}\ar[d]^-{\alpha_1} & {\shf_0}\ar[d]^-{\alpha_0}\\
_{\opb f}\shf'_1\ar[r]^-{m'} & {\shf'_0}.
}$$

\noindent  
The assignement $V\mapsto \mdso[\lienR](V)$ defines a $\opb p\kappa$-prestack 
$\mdso[\lienR]$ on $X_0$. Moreover, the position $(\shf,m)\mapsto \shf$ defines
a natural functor of prestacks $\mdso[\lienR]\to\mds[\shr]$.

\noindent
One checks that $\mds[\lienR]:=\oim{p}\mdso[\lienR]$ is a $\kappa$-stack on 
$X$, that the natural adjunction functor $F_{\lienR}\colon\opb{p}\mds[\lienR]
\to \mds[\shr]$ is an equivalence of $\opb p\kappa$-stacks on $X_0$ and that 
there exists an isomorphism  of functors 
$\alpha_{\lienR}\colon _{\opb f}(\cdot) \circ (F_{\lienR})_1 \isoto 
(F_{\lienR})_0$, these data satisfying the desired property. 
Moreover, one can show that the triplet $(\mds[\lienR],F_{\lienR},
\alpha_{\lienR})$ is unique up to equivalence.
\end{proof}

\begin{remarks}
\label{rmks}
(i) Roughly speaking, a $\kappa$-lien on $X_0$ is the data of sheaves 
of $\kappa|_{U_i}$-algebras $\shr_i$ on $U_i$ and of isomorphisms of 
$\kappa|_{U_{ij}}$-algebras $f_{ij}: \shr_j|_{U_{ij}}\isoto\shr_i|_{U_{ij}}$ 
such that $f_{ij}\circ f_{jk}=\Ad(a_{ijk})\circ f_{ik}$ on $U_{ijk}$ 
for invertible sections $a_{ijk}$ of $\shr_i|_{U_{ijk}}$. 
Hence the stack $\mds[\lienR]$ constructed in Theorem \ref{th:elienK} is a 
stack of twisted modules (see for example \cite{K-S2}), {\em i.e.,} 
$\mds[\lienR]|_{U_i}\simeq \mds[\shr_i]$.

\noindent
(ii) The section $c$ in \eqref{eq:c} defines a Cech cohomology class in 
$\check H^3 (X,\kappa^{\times})$. Indeed $c$ is a 3-cocyle, since one has the 
following chain of equalities on $\sect(X_4;\kappa^\times)$:
\eqn
&& c_{0124} c_{0234} = \opb a_{014} f_{01}(\opb{a_{124}}) a_{012} a_{024}
\opb a_{024} f_{02}(\opb a_{234}) a_{023} a_{034} \\
&& = \opb a_{014} f_{01}(\opb a_{124}) f_{01} \circ f_{12}
(\opb a_{234}) 
a_{012} a_{023} a_{034} \\
&& = \opb a_{014} f_{01}(\opb a_{124} f_{12}(\opb a_{234}) a_{123}) 
f_{01}(\opb a_{123}) a_{012} a_{023} a_{034} \\
&& = \opb a_{014} f_{01}(\opb a_{134} c_{1234}) a_{013}  c_{0123} a_{034} = 
c_{1234} c_{0134} c_{0123}.
\eneqn

\noindent
(iii) Definitions \ref{def:lien1} and \ref{def:lien2} are adapted from 
\cite{Gi} (see also \cite{Br}). In particular, the notion of an effective lien 
is a restrictive version of that of a realizable lien, {\em i.e.} a lien 
with $c$ cohomologous to 1.

\noindent
(iv) For simplicity of the exposition, the definition of a lien given 
here deals with a fixed open covering $\shu=\{U_i\}_{i\in I}$ of 
the topological space $X$. Let $\shu'=\{U'_j\}_{j\in J}$ be a finer open 
covering and define the corresponding diagram $X'_{\bullet}\to X$. 
Let $\lienR$ be a $\kappa$-lien on $X_{\bullet}\to X$.
Hence the natural map $r_0\colon X'_0\to X_0$ induces a lien 
$\lienR '=(\opb{r_0}\shr,\opb{r_1}f, \opb{r_2} a)$ on $X'_0$, where the 
$r_i$'s are the induced maps on the $X_i$'s. 
One checks that if $\lienR$ is effective, then $\lienR'$ is effective and they 
define equivalent stacks on $X$. More precisley, the equivalence 
$\mds[\lienR]\isoto\mds[\lienR']$ is given by $(\shf,m) \mapsto (\opb{r_0}\shf,
\opb{r_1} m)$. 
\end{remarks}

\begin{definition}\label{def:iso}
Let $\lienR=(\shr,f,a)$ and $\lienS=(\shs,g,b)$ be $\kappa$-liens on 
$X_{\bullet}\to X$. 
An isomorphism of $\kappa$-liens $\lien u\colon \lienR\to \lienS$ is pair 
$\lien u=(u,l)$ where $u:\shr\isoto\shs$ is an 
isomorphism of $\opb p \kappa$-algebras on $X_0$ and $l$ is a section in 
$\sect(X_1,\shs_0^\times)$ such that on $X_1$
\eq\label{eq:isoliens}
&&g \circ u_1 = \Ad(l)\circ u_0\circ f. 
\eneq
If $\lienT=(\sht,h,e)$ is another $\kappa$-lien and 
$\lien t=(t,n)\colon \lienS\to \lienT$ is an isomorphism of 
$\kappa$-liens, the composition is defined by $\lien t\circ\lien u:=(t\circ u,
n t_0(l))$.
\end{definition}

On $X_2$ one has
\eqn
g_{01}\circ g_{12}\circ u_2  = \Ad(b)\circ g_{02}\circ u_2 = 
\Ad(bl_{02})\circ u_0 \circ f_{02}
\eneqn
and 
\eqn
&& g_{01}\circ g_{12}\circ u_2 = \\
&& g_{01}\circ \Ad(l_{12}) \circ u_1 \circ f_{12} 
= \Ad(g_{01}(l_{12}))\circ \Ad(l_{01})\circ u_0\circ f_{01}\circ f_{12} = \\
&& = \Ad(g_{01}(l_{12})l_{01})\circ \Ad( u_0(a)) \circ u_0 \circ f_{02}
\eneqn
It follows that there exists a section  $d \in \sect(X_2;\kappa^\times)$, 
such that 
\eqn
&& bl_{02}= g_{01}(l_{12})l_{01} u_0(a) \cdot d \quad \mbox{ in }
\sect(X_2;\shs_0^\times)
\eneqn

\begin{definition}
One says that an isomorphism of $\kappa$-liens $\lien{u}$ is effective 
if $d=1$.
\end{definition}

\begin{proposition}\label{prop:effectiso}
Let $\lienR=(\shr,f,a)$ and  $\lienS=(\shs,g,b)$ be effective $\kappa$-liens 
on $X_{\bullet}\to X$. To an effective isomorphism of $\kappa$-liens 
$\lien{u}=(u,l)\colon \lienR\isoto \lienS$ one associates an equivalence of 
$\kappa$-stacks $F_{\lien{u}}: \mds[\lienR] \isoto \mds[\lienS]$ and an 
isomorphism of functors $\beta_{\lien{u}}\colon F_{\lienS}\circ \opb p  F_{\lien{u}} \isoto {_{\opb u}(\cdot)}\circ F_{\lienR} $ such that
$\alpha_{\lienR}\circ l\cdot \circ (\beta_{\lien{u}})_1 = (\beta_{\lien{u}})_0
\circ \alpha_{\lienS}$ (here $( \mds[\lienR], F_{\lienR},\alpha_{\lienR})$ and 
$(\mds[\lienS], F_{\lienS},\alpha_{\lienS})$ are as in Theorem \ref{th:elienK},
and  $l\cdot$ denotes the isomorphism of functors defined by left 
multiplication by $l$). 
Moreover, the datum of $(F_{\lien{u}},\beta_{\lien{u}})$ is unique up to 
unique isomorphism.
\end{proposition}

\begin{proof}
To $(\shf,m)\in\mdso[\lienR](V)$ on associates the pair 
$(_{\opb u}\shf,m\circ l \cdot)$, where the isomorphism $m\circ l \cdot$ 
stands for the following chain of isomorphisms of $\shs_0$-modules on $V_1$
\eqn
  _{\opb g}(_{\opb u}\shf)_1 \simeq 
  {}_{\opb g}(_{\opb u_1}\shf_1) 
   \to[l \cdot] {}_{\opb u_0}(_{\opb f}\shf_1)
   \to[m] {}_{\opb u_0}\shf_0
   \simeq  (_{\opb u}\shf)_0 ,
\eneqn
where $l \cdot$ denotes left multiplication by $l$ in 
$_{\opb f}\shf_1$.
Let us checks that this position gives a well defined functor of prestacks
$(F_u)_0\colon\mdso[\lienR]\to\mdso[\lienS]$.
It is enough to note that the following diagram is commutative
$$
\xymatrix{ _{\opb g_{01}}(_{\opb g_{12}}(_{\opb u}\shf)_2)
\ar[rr]^-{b\cdot} \ar[d]^-{g_{01}(l_{12})\cdot} & & 
_{\opb g_{02}}(_{\opb u_{2}}\shf_2) \ar[d]^-{l_{02}\cdot}\\
_{\opb g_{01}}(_{\opb u_1}(_{\opb f_{12}}\shf_2)) \ar[r]^-{l_{01}\cdot} 
\ar[d]^-{m_{12}} & _{\opb u_0}(_{\opb f_{01}}(_{\opb f_{12}}\shf_2)) 
\ar[r]^-{u_0(a) \cdot} \ar[d]^-{m_{12}} &  _{\opb u_0}(_{\opb f_{02}}\shf_2) 
\ar[d]^-{m_{02}} \\
_{\opb g_{01}}(_{\opb u_1}\shf_1) \ar[r]^-{l_{01}\cdot} & 
_{\opb u_0}(_{\opb f_{01}}\shf_1) \ar[r]^-{m_{01}} & 
{(_{\opb u}\shf)_0}. }
$$
Set $F_u:=p_*(F_u)_0$. This is an equivalence of $\kappa$-stacks, since a 
quasi-inverse is given by $F_{\opb u}$, and the natural morphism of functors 
$\beta_{\lien{u}}\colon F_{\lienS}\circ \opb p  F_{\lien{u}} \to 
{_{\opb u}(\cdot)}\circ F_{\lienR} $ is an isomorphism satisfying the 
desired property. 
\end{proof}

\begin{remark}
In the commutative case, an effective $\kappa$-lien on $X_{\bullet}\to X$ is 
nothing but the datum of a 2-cocycle $a\in\sect(X_2,\kappa^\times)$. 
Then there exists an effective isomorphism between two 2-cocyles 
if and only if they are cohomologous. 
\end{remark}

\begin{remark}\label{rem:groupoid0}
Consider  a diagram 
in the category of topological spaces
\eq\label{diag2}
X_\bullet:=
\xymatrix{  
X_2 
\ar@<-1.5ex>|-{p_{02}}[rr]\ar[rr]|-{p_{12}}\ar@<1.5ex>[rr]|-{p_{01}}&&
X_1 \ar@<-1.ex>[rr]|-{p_{1}}\ar@<1.ex>[rr]|-{p_{0}}&&
X_0}
\eneq
such that
\eq\label{eq:diag2}
&&\left\{ \begin{array}{l}
p_0\circ p_{01}=p_0\circ p_{02},\\
p_1\circ p_{01}=p_0\circ p_{12},\\
p_1\circ p_{12}=p_1\circ p_{02}.
\end{array}\right.
\eneq


\noindent
Then one can extend the notions of a lien and of an isomorphism of liens 
on $X_\bullet$. Moreover, there is a natural notion of continuous map between 
diagrams $r : Y_{\bullet} \to X_{\bullet}$ as above. Then, 
if $\lienR=(\shr,f,a)$ is a lien on $X_{\bullet}$, one defines the lien 
$\opb r\lienR:=(\opb r_0 \shr,\opb r_1 f,\opb r_2 a)$ on $Y_{\bullet}$.


These considerations suggest that Theorem \ref{th:elienK} extends to the
case of Lie groupoids or better, differentiable stacks 
(see \cite{Mo2}, \cite{B-X}).

\end{remark}

\section{Symplectic and contact geometry}\label{section:geometry}

Let $X$ be a complex manifold. 
Recall that a local model of a homogeneous symplectic manifold is an open 
subset $V$ of the cotangent bundle $T^*X$ equipped with the canonical $1$-form 
$\alpha_{T^*X}$. If $x=(x_1,\dots,x_n)$ is a local coordinate system on $X$, 
we denote by $(x;\xi)=(x_1,\dots,x_n;\xi_1,\dots,\xi_n)$ the associated 
coordinates on $T^*X$. Then $\alpha_{T^*X}=\sum_{i=1}^n\xi_idx_i$. 
A local model of a contact manifold is an open subset of the projective 
cotangent bundle $P^*X=\dot{T}^*X/\C^\times$, equipped with its canonical 
$1$-form. A local model of a symplectic manifold is an open subset 
$U$ of $T^*X$ equipped with the symplectic $2$-form  
$\omega_{T^*X}=d\alpha_{T^*X}$. 

We shall use the terminology ``a symplectic isomorphism'' to denote an
isomorphism of complex symplectic manifolds, and similarly for 
``a contact isomorphism'' or ``a homogeneous symplectic isomorphism''.

We denote by $\C$ the complex line endowed with a holomorphic coordinate
$t$ and by $(t;\tau)$ the associated coordinates on $T^*\C$. Hence,
$\tau$ defines a map $\tau:\dot{T}^*\C\to\C^\times$.

The embedding $T^*X\times\dot{T}^*\C\hookrightarrow\dot{T}^*(X\times\C)$ 
defines the open embedding 
$(T^*X\times\dot{T}^*\C)/\C^\times\hookrightarrow P^*(X\times\C)$. 
Set for short:
\eqn
&&\dot{P}^*(X\times\C)=(T^*X\times\dot{T}^*\C)/\C^\times.
\eneqn

\begin{remark}
There is a natural isomorphism $\dot {P}^*(X\times\C)\simeq T^*X\times \C$
given by $(p,(t;\tau))\mapsto (p \opb \tau,t)$. Hence, one can identify 
$\dot {P}^*(X\times\C)$ with the space $J^1 X$ of 1-jets of holomorphic
functions on $X$. 
\end{remark}

We shall consider the maps:
\eq\label{diag:1}
\xymatrix{
{\C^\times}  &{T^*X\times\dot{T}^*\C}\ar[l]_-\tau \ar[r]^-{\gamma}
       &{\dot{P}^*(X\times\C)} \ar[r]^-{\rho}&T^*X }
\eneq
where $\rho$ is the map $(p,(t;\tau))\mapsto p \opb \tau$. 

\begin{definition}
Consider a contact isomorphism $\psi:V_X\isoto V_Y$, where $V_X$ \lp
resp. $V_Y$\rp\, is an open subset of 
$\dot{P}^*(X\times \C)$ \lp resp. of $\dot{P}^*(Y\times \C)$\rp.
We say that $\psi$ is a $\tau$-preserving contact isomorphism,
if it lifts as a homogeneous symplectic isomorphism 
$\widetilde{\psi}$ making the diagram below commutative:
\eqn
\xymatrix{
{\dot{P}^*(X\times \C)\supset V_X}\ar[rr]^-\psi
                                      &&{V_Y\subset\dot{P}^*(Y\times\C)}\\
{T^*X\times \dot{T}^*\C\supset \opb{\gamma}(V_X)}\ar[u]_\gamma
                         \ar[rr]^-{\widetilde{\psi}}\ar[rd]_-\tau
 &&{\opb{\gamma}(V_Y)\subset T^*Y\times\dot{T}^*\C}\ar[u]_\gamma\ar[ld]^-\tau\\
             &{\C^\times}&  }
\eneqn
\end{definition}

\begin{lemma}\label{lemma:geomconttrans}
Let $U_X$ \lp resp. $U_Y$\rp\, be an open subset of $T^*X$ \lp
resp. $T^*Y$\rp\, and let 
$\phi:U_X\isoto U_Y$ be a symplectic isomorphism. 
Then, locally on $U_X$, there exists a $\tau$-preserving contact isomorphism 
$\psi:{\opb{\rho}(U_X)}\to {\opb{\rho}(U_Y)}$ 
making the diagram below commutative 
\eq\label{diag:2}
\xymatrix{
{T^*X\supset U_X}\ar[r]^-\phi&{U_Y\subset T^*Y}\\
{\dot{P}^*(X\times \C)\supset \opb{\rho}(U_X)}\ar[u]_-\rho\ar[r]^-\psi
           &{\opb{\rho}(U_Y)\subset\dot{P}^*(Y\times\C)}\ar[u]_-\rho  }
\eneq
\end{lemma}

\noindent
One shall be aware that $\psi$ commutes with $\tau$, but not with $t$
in general.

\begin{proof}
Let us denote by $(x;u)$ a local symplectic coordinate system on $U_X$
and let $(y;v)=\phi(x;u)$. Then $du\wedge dx=dv\wedge dy$, hence 
$d(udx)=d(vdy)$ (we write for short $udx$ instead of $\langle u,dx\rangle$). It
follows that
$udx=vdy+ da(y,v)$ for some locally defined function $a(y,v)$. 

\noindent
Set $\xi=\tau u, \eta=\tau v$ and consider the map $\psi$ such that
$$(y,s;\eta,\tau)=\psi(x,t;\xi,\tau)$$
for an $s$ to be calculated below. We have
 \eqn 
\tau dt+\xi dx=\tau(dt +udx)=\tau(dt+vdy+da)=\tau d(t+a)+\eta dy.
\eneqn
Hence, choosing $s=t+a(y,\eta \opb\tau)$, the map $\psi$ is a 
$\tau$-preserving contact isomorphism.
\end{proof}

Note that, if $\phi$ is an homogeneous symplectic isomorphism, then the 
function $a$ constructed above is globally defined and locally constant on 
$U_X$. In particular, if $\phi$ is the identity on $U_X$, then 
a $\tau$-preserving contact isomorphism $\psi$ is nothing but a translation on 
the $t$-variable by a locally constant function.

\section{$\she$-modules}\label{section:Emod}

On $T^*X$ we consider the sheaves $\she_X$ and $\shhe_X$ of
microdifferential operators and formal microdifferential operators, 
respectively.
We refer the reader to \cite{S-K-K} for more details
(see also \cite{S} or \cite{K1} for an exposition).

These sheaves are constant on the fibers of the projection
$\dot{T}^*X\to P^*X$, and we shall keep the same notation to denote
their direct images on $P^*X$. Hence we regard them as 
sheaves on $\dot{T}^*X$ as well as sheaves on $P^*X$. 

The sheaf  $\she_X$ is filtered over $\Z$, and one denotes by
$\she_X(m)$ the sheaf of operators of order less than or equal to $m$. We
denote by $\sigma_m(\cdot):\she_X(m)\to \she_X(m)/\she_X(m-1)$ the
symbol map. Recall that $\she_X(m)/\she_X(m-1)\simeq \sho_{T^*X}(m)$,
the sheaf of holomorphic functions homogeneous of order $m$ 
in the fiber variable. A section $f$ of $\sho_{T^*X}(m)$ is a
holomorphic function solution of the differential equation 
$(\sum_jx_j\partial_{x_j}-m)f=0$.
Hence
\eqn
&&gr(\she_X)\simeq\bigoplus_{j\in\Z}\sho_{T^*X}(j).
\eneqn
The same result holds for $\widehat\she_X$.

In a local coordinate system $x$ on $X$, with associated coordinates
$(x;\xi)$ on $T^*X$, a formal microdifferential operator $P$ of order $m$
({\em i.e.,} a section of $\shhe_X(m)$) defined on an open subset $V$ of 
$T^*X$  has a total symbol $\sigma(P)$:
\eq\label{eq:totsymb}
&&\sigma(P)=\sum_{j=-\infty}^m p_j(x,\xi), \,\,\,p_j\in\sho_{T^*X}(j)(V).
\eneq

\begin{notation}
In a local coordinate system $x=(x_1,\dots,x_n)$ on $X$, one denotes by
$\partial_{x_j}$ (or else, $D_{x_j}$) the microdifferential operator 
with total symbol
$\xi_j$. Hence, $P$ is written as 
\eq\label{eq:totsymbbis}
&&P=\sum_{j=-\infty}^m p_j(x,\partial_x).
\eneq
\end{notation}

The product structure on $\shhe_X$ is then given by the
Leibniz formula. If $Q$ is a formal microdifferential operator of 
 total symbol $\sigma(Q)$, then 
$$
\sigma(P\circ Q)=
\sum_{\alpha\in \N^n} \frac{1}{\alpha !} 
\partial^{\alpha}_\xi\sigma(P)
\partial^{\alpha}_x\sigma(Q).
$$

\noindent
In particular, a section $P$ in $\shhe_X$ is invertible on an open subset $V$
of $T^*X$ if and only if its principal symbol is nowhere vanishing on $V$.

The ring $\she_X$ is the subring of $\shhe_X$ of operators
whose total symbol  \eqref{eq:totsymb} satisfies the estimates
\eq\label{eq:estmicrod}
&&\left\{ \begin{array}{l}
\mbox{for any compact subset $K$ of $V$ there exists a constant}\\
\mbox{$C_K>0$ such that for all $j<0$, }
\sup\limits_{K}\vert p_{j}\vert \leq C_K^{-j}(-j)!.
\end{array}\right.
\eneq

Let $\lambda$ be a complex number. Replacing $\sho_{T^*X}(j)$ with 
$\sho_{T^*X}(\lambda +j)$,  
one defines by the same procedure the sheaves $\she_X(\lambda)$ and 
$\shhe_X(\lambda)$ on $T^*X$ of operators of order $m+\lambda$. 
Note that in the $P^*X$ case, one obtains twisted sheaves.
 
From now on, we shall concentrate our study on $\she_X$, but all
results extend unchanged to $\shhe_X$.

A volume form on $X$ defines an
anti-automorphism $*:\she_X \to \oim{a}\she_X$ (recall that $a$ is the 
antipodal map on $T^*X$).
This leads to consider the sheaf of rings:
\eq\label{eq:shedens}
&&\shed[X]:=
\Omega_X^{\tens 1/2}\tens_{\sho_X}\she_X\tens_{\sho_X}\Omega_X^{\tens -1/2}.
\eneq
(Here, we write $\sho_{X}$ instead of
$\opb{\pi}\sho_{X}$ and similarly for $\Omega_X$.) Note that 
$\Omega_X^{\tens 1/2}$ and $\Omega_X^{\tens -1/2}$ are not globally
defined as sheaves but are globally 
defined as twisted sheaves. On the other-hand
$\shed[X]$ is a well-defined sheaf of rings on $T^*X$, locally isomorphic to
$\she_X$.
This sheaf on $P^*X$ (not on $T^*X$) is thus endowed with an
anti-automorphism $*$ such that $**=\id$.

\subsection*{Quantized contact transformations}
Let us briefly recall the constructions of ``quantized contact
transformations'' of \cite{S-K-K}.

Assume to be given open subsets 
$V_X$ of $T^*X$, $V_Y$ of $T^*Y$\ and a homogeneous 
symplectic isomorphism $\psi:V_X\isoto V_Y$. Let
$\Lambda\subset V_X\times V_Y^a$ denote the Lagrangian submanifold
associated with $\psi$, that is, the image of the graph of $\psi$ by
the antipodal map $a:V_Y\to V_Y^a$. Locally on $\Lambda$, there exists
a left ideal $\shi$ of $\she_{X\times Y}$ such that its symbol ideal is
reduced and coincides with the defining ideal of $\Lambda$ in 
$\bigoplus_{j\in\Z}\sho_{T^*X}(j)$.

Then for each $P\in\she_X$ (defined in a neighborhood of $p\in V_X$) 
there exists a unique $Q\in\she_Y$ 
(defined in a neighborhood of $\psi(p)\in V_Y$) such that $(P-Q)\in\shi$. 
The correspondance $P\mapsto Q$ is an anti-isomorphism of
$\C$-algebras. Composing it with the anti-isomorphism $Q\mapsto Q^*$
associated with a volume form on $Y$, we find an isomorphism of
$\C$-algebras
\eq\label{eq:qct1}
&&\Psi:\oim{\psi}(\she_X|_{V_X})\isoto\she_Y|_{V_Y}.
\eneq

The same construction holds on projective cotangent bundles. It also
holds replacing $\she$ with $\shde$.

\begin{definition}
\bnum
\item
Let $\psi$ be a homogeneous symplectic isomorphism. The isomorphism $\Psi$ 
in \eqref{eq:qct1} is called a homogeneous
quantized contact transformation above $\psi$ (a homogeneous 
QCT, for short).
\item
If $\psi$ is a contact transformation ({\em i.e.} $V_X$ and $V_Y$ are open 
subsets of projective cotangent bundles),
one calls $\Psi$ a quantized contact 
transformation (a QCT, for short).
\item
We keep the same terminology when $\she$ is replaced with 
$\shde$.
\enum
\end{definition}

\begin{lemma} {\rm \cite{S-K-K}, \cite{K1}, \cite{K2}.}\label{le:adpdc}
\bnum
\item 
Let $\Psi$ be a QCT above the identity. 
Locally there exists $\lambda\in\C$ and an invertible operator
$P\in\she_X(\lambda)$ such that 
$\Psi=\Ad(P)$. 
\item
Let $P\in\she_X$ be of order $\lambda$ and invertible. 
There exists a unique (up to sign) 
invertible operator $Q$ of order $\frac{1}{2}\lambda$ 
such that $P=Q\circ Q$. If $P\in\shde_X$, $P$ is of order $0$ and $P=P^*$, 
then $Q=Q^*$.
\item
Let $\psi$ be a contact isomorphism. Locally 
there exists a QCT
\eq
&&\Psi:\oim{\psi}(\shed[X]\vert_{V_X})\isoto \shed[Y]\vert_{V_Y} 
\eneq
commuting with $*$.
\enum
\end{lemma}
\begin{proof}
(i) We shall not give the proof here and refer to  \cite{K1}.

\noindent
(ii) (a) Unicity. First assume $P$ has order $0$. Let $P=Q^2=Q_0^2$.
Then $\sigma_0(Q_0)=\pm\sigma_0(Q)$ and we may assume there is
equality. Set $Q_0=Q+R$ and let $m$ denote the order of $R$. 
One has $Q_0^2=(Q+R)^2$ and therefore $QR+RQ+R^2=0$. 
Since $2m<m$, it follows that $\sigma_m(R)=0$, hence $R=0$.

\noindent
The general case follows by adding a dummy  variable $\partial_t$,
and applying the preceding result to $\partial_t^{-\lambda}P$.

\noindent
(ii) (b) Existence. First assume $P$ has order $0$. 
Let $p_0\in T^*X$ and assume 
$\sigma_0(P)(p)\neq 0$. Consider a dummy
variable $t$ and set $t_0=(\sigma_0(P)(p_0))^{\frac{1}{2}}$. 
Since $t^2-\sigma_0(P)$ has a simple root at
$(t_0,p_0)$, the Weierstrass Preparation Theorem for
microdifferential operators of $\cite{S-K-K}$ (see also \cite{S} Chapter I 
\S 2 for an exposition) allows us to write uniquely
\eqn
&& t^2-P= G(t-R)
\eneqn
with $G$ invertible at $(t_0,p_0)$ and $R$ not depending on $t$. 
Therefore $P-R^2=(t-R)(t+R-G)$, and this operator not depending on
$t$, we find $G=t+R$. Hence, $P=R^2$.

\noindent
The general case follows by adding a dummy  variable $\partial_t$.
Then $\partial_t^{-\lambda} P=R^2$. Set 
$Q=\partial_t^{\frac{\lambda}{2}}R$. Then $P=Q^2$, and the unicity of
this decomposition shows that $Q$ does not depends on $\partial_t$. 

\noindent
(ii) (c) Since $P=P^*$, we get $Q^*Q^*=QQ$, and by the unicity of the
decomposition, $Q=\pm Q^*$. Since $\sigma_0(Q)=\sigma_0(Q^*)$, we get
$Q=Q^*$.

\noindent
(iii) First choose a QCT $\Psi^{\dagger}$. 
Set $\Psi^{\ddagger}=\Psi^{\dagger -1}\circ *\circ\Psi^{\dagger}\circ *$. 
By (i), $\Psi^{\ddagger}=\Ad(R)$ for an invertible operator
$R$ of order $\lambda$. Hence
\eqn
 *\circ\Psi^{\dagger}\circ * & = &\Psi^{\dagger}\circ \Ad(R),\\
 \Psi^{\dagger} & = & *\circ\Psi^{\dagger}\circ *\circ \Ad(\opb{{R^*}})\\
                & = &\Psi^{\dagger}\circ \Ad(R)\circ \Ad(\opb{{R^*}})\\
                & = &\Psi^{\dagger}\circ\Ad(R \circ\opb{{R^*}}).
\eneqn
Therefore, $R=cR^*$ for some non-zero constant $c$, that we may assume 
equal to $1$. Hence using (ii), we may write
$R=Q\circ Q^*$ and we get
$\Psi^{\ddagger}=\Ad(Q)\Ad(Q^*)$. Let $\Psi=\Psi^{\dagger}\circ \Ad(Q)$. 
Then, one has
$$
*\circ\Psi^{\dagger} \circ * = \Psi \circ \Ad(Q^*).
$$
Since $\Ad(Q)\circ * = *\circ \Ad(\opb{{Q^*}})$,
one gets
$$
*\circ \Psi \circ * = \Psi.
$$
\end{proof}

Denote by $\sha ut_* (\shed[X])$ the sheaf of QCT above the identity
commuting with $*$ and set 
\eqn
&&(\shed[X])_*=\{P\in\shed[X];P\text{ has order }0, \sigma_0(P)=1, PP^*=1\}
\subset (\shed[X])^{\times}.
\eneqn 

\begin{lemma}{\rm \cite{K2}}\label{le:auts}
The morphism $\Ad$ induces an isomorphism of groups on $T^*X$ 
\eqn
(\shed[X])_* \isoto[\Ad] \sha ut_*(\shed[X]).
\eneqn
\end{lemma}

\begin{proof}
Let $\Psi$ be a QCT above the identity commuting with $*$. By Lemma 
\ref{le:adpdc}, locally there exists an invertible $P\in \shed[X]$
of order $\lambda$ such that $\Psi=\Ad(P)$. One has 
\eqn
&& *\circ \Ad(P)\circ *=\Ad(\opb{{P^*}}).
\eneqn
Since $*\circ \Psi=\Psi\circ *$, we get $\Ad(\opb{{P^*}})=\Ad(P)$, 
and $C := P^*P$ is an invertible element in the center 
of $\shed[X]$. Therefore, $P$ has order $0$ and $\sigma_0(P)$ is 
a non-zero constant. Then  we may suppose $P$ of principal symbol $1$.
\end{proof}

\begin{remark}
The results of this section still hold when replacing $\she$ with $\shhe$.
\end{remark}

\section{Quantization of complex contact manifolds}\label{section:Cquant}

In this section, we recall Kashiwara's theorem.

\begin{theorem}{\rm (M. Kashiwara \cite{K2}.)}\label{th:qccm}
Let $\sty$ be a complex contact manifold.
There exists canonically 
a $\C$-abelian stack $\mds[\shde,\sty]$ on $\sty$ such 
that if $V\subset\sty$ is an open subset  isomorphic by a contact
transformation $\psi$ to an open subset $V_X \subset P^*X$, then 
$\mds[\shde,\sty]\vert_V$ is equivalent by $\psi$ to the stack 
$\mds[{\shed[X]\vert_{V_X}}]$. 
\end{theorem}

\begin{definition}
We call $\mds[\shde,\sty]$ the stack of microdifferential modules 
on $\sty$. 
\end{definition}

\begin{proof}
There exists an open covering $\shv=\{V_i\}_{i\in I}$ of $\sty$ and for each
$i\in I$, a contact open embedding $\psi_i:V_i\hookrightarrow P^*X_i$ for 
some projective cotangent bundle $P^*X_i$. 
Set $\shed[V_i]=\opb{\psi_i}\shed[X_i]$. Then $\shed[V_i]$ is a sheaf of
$\C$-algebras on $V_i$ endowed with a filtration and an anti-involution $*$.

\noindent
Consider the contact isomorphism $\psi_{ij}=\psi_i\circ \opb \psi_j : 
\psi_j(V_{ij})\isoto \psi_i(V_{ij})$. 
After shrinking the covering $\shv$, we may 
assume by Lemma \ref{le:adpdc} that there exist QCT's above the $\psi_{ij}$'s 
and hence isomorphisms of $\C$-algebras
\eq\label{gluedatasymp1}
&&\Psi_{ij}:\shed[V_j]\vert_{V_{ij}}\isoto\shed[V_i]\vert_{V_{ij}},
\eneq
these isomorphisms commuting with $*$. 

\noindent
Now we follow the notations of Section \ref{section:stacks} for $X=\sty$.
Set $\sty_0=\bigsqcup_iV_i$,  $\sty_1=\sty_0\times_{\sty}\sty_0$, etc. Let 
$j_{V_i}:V_i\hookrightarrow \sty_0$ be the natural map. Set 
\eqn
&&\shed[]=\oplus_{i\in I}\eim{{j_{V_i}}}\shed[V_i].
\eneqn
Then $\shed[]$ is a sheaf of central 
$\C$-algebras on $\sty_0$ endowed with an anti-involution $*$. 
The $\Psi_{ij}$'s induce a $\C$-algebra isomorphism 
$\Psi:\shed[1]\isoto\shed[0]$ commuting $*$. Hence, 
by Lemma \ref{le:auts}, after shrinking again the covering $\shv$ 
there exists a unique $P\in\sect(\sty_2;\shed[0])$ 
of order $0$ such that $\sigma_0(P)=1$, $PP^*=1$ and  
\eq\label{cocyle5}
\Psi_{01}\circ \Psi_{12}=\Ad(P)\circ \Psi_{02}.
\eneq

\noindent
Since $P$ is unique, $\lien{E}:=(\shed[],\Psi,P)$ is an effective $\C$-lien
on $\sty_{\bullet}\to\sty$, and it remains to apply Theorem \ref{th:elienK} to
get an abelian $\C$-stack $\mds[\lien{E}]$ on $\sty$.

Let $\shv'=\{V'_j\}_{j\in J}$ be another open covering of $\sty$ and 
$\psi'_j:V'_j\hookrightarrow P^*X'_j$ be a contact open embedding for each 
$j\in J$. By Remark \ref{rmks} (iii), it is not restrictive to assume 
$\shv = \shv'$. Setting ${\shed[V_i]}'=\opb{{\psi'_i}}\shed[X'_i]$ 
and proceeding 
as above, we get another effective $\C$-lien $\lien{E}':=({\shed[]}',\Psi',P')$
on $\sty_{\bullet}\to\sty$. 
Consider the contact isomorphism $\psi'_i \circ \opb {\psi_i} : 
\psi_i(V_i)\isoto \psi'_i(V_i)$. 
After shrinking the covering $\shv$, we may 
assume by Lemma \ref{le:adpdc} that there exist QCT's above the 
$\psi'_i \circ \opb {\psi_i}$'s and hence an isomorphism of $\C$-algebras
\eq
&&\Upsilon:\shed[]\isoto{\shde}'
\eneq
commuting with $*$.
By Lemma \ref{le:auts}, after shrinking again the covering $\shv$, 
there exists a unique $Q\in\sect(\sty_1;{\shde}')$ 
of order $0$ such that $\sigma_0(Q)=1$, $QQ^*=1$ and  
\eq
\Psi' \circ \Upsilon_1 = \Ad(Q) \circ\Upsilon_0\circ \Psi.
\eneq
Since $Q$ is unique, the pair $(\Upsilon,Q)$ defines an effective isomorphism 
of $\C$-liens $\lien{E}\isoto \lien{E}'$. By Proposition
\ref{prop:effectiso}, the stacks $\mds[\lien{E}]$ and $\mds[\lien{E}']$ are 
equivalent.

\noindent
Hence the stack above constructed depends only on $\sty$ and on the algebra
$\shed[]$, and it makes sense to denote it by $\mds[\shde,\sty]$. 
\end{proof}

\begin{remark}\label{rem:ehat}
The results of this section still hold when replacing $\she$ with $\shhe$.
\end{remark}

\section{Quantization of involutive submanifolds}\label{section:Invquant}

We keep the notations of \S \ref{section:Cquant} and
consider a complex contact manifold $\sty$.
In this section, we consider  a smooth 
regular involutive submanifold $\Lambda$ 
and we denote by 
\eq
&&\iota\colon\Lambda \hookrightarrow \sty
\eneq 
the inclusion morphism. Recall that one says that $\Lambda$ is involutive 
if for any pair of holomorphic functions $(f,g)$ vanishing on $\Lambda$, 
their Poisson bracket $\{f,g\}$ vanishes on $\Lambda$ and recall that 
$\Lambda$ is regular involutive if moreover the canonical
1-form on $\sty$  does not vanish on $\Lambda$. 

It is well-known that locally, two smooth regular involutive
manifolds of the same codimension may be interchanged by a complex contact
isomorphism. In particular, locally on $\Lambda$, we may assume that
\eq\label{eq:lambdalocal}
\sty=P^*X, &X=Y\times Z,&\Lambda=P^*Y\times Z. 
\eneq

Let $b$ denotes the bicharacteristic relation on $\Lambda$, which identifies a 
bicharacteristic leaf to a point. We denote by 
\eq
&&\beta:\Lambda\to\Lambda/b
\eneq 
the projection, and if $U$ is open in $\Lambda$ we denote by 
$\beta_U:U\to U/b$ 
the restriction of $\beta$ to $U$. Note that for $U$ small enough,
 $U/b$ has the 
structure of complex contact manifold. 

We can now formulate a variant of Theorem \ref{th:qccm}.

\begin{proposition}\label{pr:qccmB}
Let $\sty$ be a complex contact manifold and let $\Lambda$ be 
a smooth regular involutive submanifold of $\sty$.
There exists canonically a $\C$-abelian stack 
$\mds[\opb{\beta}\shde,\Lambda]$ on $\Lambda$ such that 
if $U\subset\Lambda$ is an open subset  and $U/b$ is a contact manifold 
isomorphic by a contact transformation $\psi$ to an open subset $V\subset P^*Y$,
then $\mds[\opb{\beta}\shde,\Lambda]\vert_U$ is equivalent to 
the stack $\mds[{\opb{\beta_U}\opb{\psi}\shed[Y]\vert_V}]$. 
\end{proposition}

\begin{proof}
Consider an open covering $\shu=\{U_i\}_{i\in I}$ of $\Lambda$
such that for each $i\in I$, $U_i/b$ is a smooth
complex contact manifold and there exists a contact embedding
$\psi_i:U_i/b\hookrightarrow P^*Y_i$
for some projective cotangent bundle  $P^*Y_i$.
Set $V_i=U_i/b$, $V_{ij}=(U_i\cap U_j)/b$, $V_{ijk}=(U_i\cap U_j\cap U_k)/b$,
$\sty_0=\bigsqcup V_i$, $\sty_1=\bigsqcup V_{ij}$, 
$\sty_2=\bigsqcup V_{ijk}$.
We find a diagram $\sty_\bullet$ as in \eqref{diag2} and
after shrinking the covering $\shu$, 
we find an effective $\C$-lien $\lien{E}:=(\shed[],\Phi,P)$ on $\sty_\bullet$.

\noindent
Set $\Lambda_0=\bigsqcup U_i$, $\Lambda_1=\bigsqcup U_{ij}$, 
$\Lambda_2=\bigsqcup U_{ijk}$, and define the diagram 
$\Lambda_\bullet\to\Lambda$ similarly. The projection $\beta$ defines a 
continuous map of diagrams $\beta : \Lambda_{\bullet}\to\sty_{\bullet}$.
Hence $\opb\beta \lien{E}$ is an effective $\C$-lien on 
$\Lambda_\bullet\to\Lambda$, and it remains to apply Theorem \ref{th:elienK}.

\noindent
The proof of the canonicity of the stack $\mds[\opb{\beta}\shde,\Lambda]$ goes 
along the same lines as in the proof of Theorem \ref{th:qccm}.
\end{proof}

Consider the situation of \eqref{eq:lambdalocal}: 
$\sty=P^*X$, $X=Y\times Z$, $\Lambda=P^*Y\times Z$ and 
denote by $\beta:P^*Y\times Z\to P^*Y$ the projection. 
Let $\shl_{\Lambda}:=\she_Y\underline{\etens} \sho_Z$ 
where the external product in taken in the category of
$\she$-modules. Such an $\she_X$-module is called  
simple along $\Lambda$. Then
$\she nd_{\she_X}(\shl_\Lambda)\simeq \opb{\beta}\she_Y$ and 
 $\shl_{\Lambda}$ is faithfully flat over 
$\opb{\beta}\she_Y$. This suggests another method to quantize the
involutive submanifold $\Lambda$, namely by glueing the 
categories of modules over $\she nd_{\she_X}(\shl_\Lambda)$, for local
choices of simple modules $\shl_{\Lambda}$.

\begin{remark}\label{rem:ehat2}
The results of this section still hold when replacing $\she$ with $\shhe$.
\end{remark}

\section{$\shw$-modules} \label{section:Wmod}

Let $X$ be a complex manifold, and let $\C$ denote as above the
complex line endowed with the holomorphic coordinate $t$.
Recall the map $\rho$ constructed in \eqref{diag:1}
\eqn
&&\rho:\dot{P}^*(X\times\C)\to T^*X.
\eneqn

\begin{definition}
\bnum
\item
On $\dot{P}^*(X\times\C)$,
we denote by $\shet[X]$ the subsheaf of rings of $\she_{X\times\C}$ 
consisting of microdifferential operators $P$ such that 
\eq
&&[P,D_t]=0.
\eneq
For $m\in\Z$, we set $\shet[X](m)=\shet[X]\cap\she_{X\times\C}(m)$.
\item
On $T^*X$ we define the sheaf of rings
\eqn
&&\shw_X=\oim{\rho}\shet[X].
\eneqn
We set $\shw_X(m)=\oim{\rho}\shet[X](m)$.
\item
Replacing $\she$ with  $\shde$, $\shhe$ and $\shhed[]$, 
one constructs similarly the sheaves of rings $\shwd[X]$,
$\shhw_X$ and $\shhwd[X]$ on $T^*X$.
\enum
\end{definition}

After choosing a local coordinate system $x$ on $X$, the 
microdifferential  operator
$P(x,t;D_x,D_t)\in \shet[X]$ does not depend on $t$ and its total symbol 
may be written as a serie
\eq
&&\sigma(P)=\sum_{j= -\infty}^m p_j(x;\xi,\tau)
\eneq
where the $p_j$' are holomorphic functions defined on some conic open
subset $V$ of $T^*(X\times \C)$, homogeneous of order $j$ with respect
to $(\xi,\tau)$.  
Setting $u=\xi \opb\tau$, $\tilde p_j(x;u)= p_j(x;\xi \opb\tau,1)$, we get that
a section $P$ of $\shw_X$ on an open subset $U$ of $T^*X$ 
has a total symbol
\eq\label{eq:totsymbsymb}
&&\sigma(P)=\sum_{j=-\infty}^m \tilde p_j(x;u)\tau^{j}
\eneq
where the $\tilde p_j$'s are holomorphic (but no more homogeneous) on 
$U$.
These functions should satisfy the follow estimates
\eq\label{eq:estmicrod}
&&\left\{ \begin{array}{l}
\mbox{for any compact subset $K$ of $U$ there exists a constant}\\
\mbox{$C_K>0$ such that for all $j<0$, }
\sup\limits_{K}\vert \tilde p_j\vert \leq C_K^{-j}(-j)!.
\end{array}\right.
\eneq

\begin{notation}
In a local coordinate system $x=(x_1,\dots,x_n)$ on $X$, one denotes by
$\opb{\tau}\partial_{x_j}$ (or else, $\opb{\tau}D_{x_j}$) 
the operator with total symbol
$u_j$. Hence, an operator $P$ is written as 
\eq\label{eq:totsymbsymbis}
&&P=\sum_{j=-\infty}^m p_j(x,\opb{\tau}\partial_x)\tau^{j}.
\eneq
\end{notation}

The ring $\shw_X$ is filtered and 
\eq\label{eq:grsymbwo}
&& gr(\shw_X)\simeq\sho_{T^*X}[\tau,\opb\tau].
\eneq
As in the case of
microdifferential operators, the symbol of order $m$ of $P$ is denoted
$\sigma_m(P)$. This function does not depend on the local coordinate
system on $X$.
If $\sigma_m(P)$ is not
identically zero, then one says that $P$ has order $m$ and
$\sigma_m(P)$
is called the principal symbol of $P$. 

The product in $\shw_X$ (and in $\shhw_X$) is given by the Leibniz formula 
not involving the $\tau$-derivatives.
If $Q$ is an operator of total symbol $\sigma(Q)$, then 
$$\sigma(P\circ Q)=\sum_{\alpha\in\N^n} \dfrac{\tau^{-\vert\alpha\vert}}
{\alpha !} \partial^{\alpha}_u\sigma(P)\partial^{\alpha}_x\sigma(Q).
$$ 
In particular, a section $P$ in $\shw_X$ is invertible on an open subset $U$ 
of $T^*X$ if and only if its principal symbol is nowhere vanishing on $U$.

\begin{remarks}
(i) The rings $\shw_{X}$ and $\shhw_{X}$ have already been introduced (when 
$X=\C$) in \cite{AKKT}. They are denoted $\she_{WKB}$ and $\shhe_{WKB}$
by these authors who call their sections, WKB-differential
operators. We shall keep this last terminology.

\noindent
(ii)
The ring $\shhw_X$ is a  semi-classical star-algebra in the 
sense of \cite{Bo}: it is locally isomorphic to 
$\sho_{T^*X}[\tau,\opb\tau]\!]$ as a $\C_{T^*X}$-module (via the total symbol)
and it is equipped with an unital associative product (the Leibniz rule) which
induces a star-product on $\sho_{T^*X}[\tau,\opb\tau]\!]$ not involving 
the $\tau$-derivatives.
\end{remarks}

\begin{definition}\label{def:cor}
\bnum
\item
One denotes by $\hcor$ the field
$\C[\tau,\opb\tau]\!]$, that is, the field of formal series  
$\sum_{j\in\Z}a_j{\tau}^j$ with $a_j=0$ for $j>\!>0$.
\item
One denotes by $\cor$ the subfield of $\hcor$ consisting 
of series $\sum_j a_j{\tau}^j$ which satisfies the estimate:
\eq\label{eq:cor}
&&\left\{ \begin{array}{l}
\mbox{there exists a constant $C>0$ such that }\\
\mbox{for all $j<0$, }\vert a_{j}\vert \leq C^{-j}(-j)!.
\end{array}\right.
\eneq
\enum
In other words, $\cor=\shw_{\rmpt}$ and $\hcor=\shhw_{\rmpt}$.
\end{definition}

The following result follows immediately from the fact that the center
of $\she_X$ and $\shhe_X$ is $\C_{T^*X}$.

\begin{lemma}\label{lemma:center}
\bnum
\item
The center of $\shw_X$ is the constant sheaf ${\cor}_{T^*X}$. 
\item
The center of $\shhw_X$ is the constant sheaf ${\hcor}_{T^*X}$.
\enum
\end{lemma}

\begin{remark}
Replacing the sheaf of rings $\she$ with the sheaf of rings $\shie$ of 
infinite order microdifferential operators, one constructs similarly the 
sheaf of rings $\shiw_X$ on $T^*X$ of infinite order WKB-differential 
operators.
\end{remark}

\subsection*{Quantized symplectic transformations}

\begin{lemma}\label{lemma:tauQCT2}
Assume to be given an open subset 
$U_X$ \lp resp. $U_Y$\rp\, of $T^*X$ \lp resp. $T^*Y$\rp\,, a symplectic 
isomorphism $\phi:U_X\isoto U_Y$ 
and a $\tau$-preserving contact isomorphism 
$\psi: \opb{\rho}(U_X)\to\opb{\rho}(U_Y)$ 
making the diagram \eqref{diag:2} commutative.
Denote by $x=(x_1,\dots,x_n)$ a local coordinate system on $X$,
by $(x,t;\xi,\tau)$ the associated local homogeneous coordinate system on
$\opb{\rho}(U_X)$ and by $(y,s;\eta,\tau)$ its image in 
$\opb{\rho}(U_Y)$ by the $\tau$-preserving contact isomorphism $\psi$
{\rm (}hence, $y_j = f_j(x,\xi,\tau)$, $\eta_j = g_j(x,\xi,\tau)$ and 
$s=t+a(x,\xi,\tau)$ {\rm )}.
Then, there locally exists a QCT above $\psi$ 
\eq\label{eq:QCST1}
&&\Psi:\oim\psi(\shed[X\times\C]\vert_{\opb{\rho}(U_X)})
    \isoto \shed[Y\times\C]\vert_{\opb{\rho}(U_Y)}
\eneq
satisfying:
\eqn
&&\Psi\mbox{ commutes with the anti-involution $*$,}\\
&&\Psi(\partial_t)=\partial_t,\quad 
\Psi(t)=t+S(x,\partial_x,\partial_t),\quad\sigma_0(S)=a ,\\ 
&&\Psi(x_j)=P_j(x,\partial_x,\partial_t),\quad\sigma_0(P_j)=f_j,\\
&&\Psi(\partial_{x_j})=Q_j(x,\partial_x,\partial_t),\quad\sigma_1(Q_j)=g_j.
\eneqn
\end{lemma}

\begin{proof}
Quantizing the contact transformation $\psi$ means finding microdifferential
operators $P_j(x,t,\partial_x,\partial_t)$ and
$S(x,t,\partial_x,\partial_t)$ of order $0$, 
$Q_j(x,t,\partial_x,\partial_t)$ and 
$T(x,t,\partial_x,\partial_t)$ of order $1$, satisfying:
\eqn
&&[P_i,Q_j]=-{\delta}_{ij},\quad [P_i,P_j]=0, \quad [Q_i,Q_j]=0\\
&&[P_i,S]=0,\quad [P_i,T]=0,\quad [Q_i,T]=0,\quad [S,T]=-1,\\
&&\sigma_0(P_j)=f_j, \quad \sigma_1(Q_j)=g_j, \quad \sigma_1(T)=\tau . 
\eneqn
One may choose $T=\partial_t$ and it follows that 
$P_j,Q_j, S-t$ will be independant of $t$ and $\sigma_0(S)=a$.
\end{proof}

Recall that $\rho:\dot{P}^*(X\times \C)\to T^*X$ is given by 
$\rho(p,(t,\tau))=(p\opb\tau)$.
Taking the direct image by $\rho$, the isomorphism $\Psi$ induces an 
isomorphism
\eq\label{eq:psi1}
&&\Phi:\oim{\phi}\oim{\rho}(\shed[X\times \C]\vert_{\opb{\rho}(U_X)})
\isoto \oim{\rho}(\shed[Y\times \C]\vert_{\opb{\rho}(U_Y)}).
\eneq

We identify $\shw_X$ with the subsheaf of algebras of 
$\oim{\rho}\she_{X\times \C}$ consisting of sections commuting with
$\partial_t$. Since $\partial_t$ is central in $\shw_X$, 
we denote it by $\tau$,
that is, we identify the operator and its symbol.
Since $\Psi(\partial_t)=\partial_t$, the isomorphism $\Phi$ induces an
isomorphism of filtered $\cor$-algebras (which we denote by the same symbol):
\eq\label{eq:QST1}
&&\Phi:\oim{\phi}(\shwd[X]\vert_{U_X})\isoto \shwd[Y]\vert_{U_Y}. 
\eneq

\begin{definition}
We call the isomorphism
$\Phi$ constructed in \eqref{eq:QST1} a quantized symplectic transformation
(a QST, for short) above $\phi$.
\end{definition}

Note that $\Phi$ depends on the $\tau$-preserving contact isomorphism $\psi$.

\begin{definition}
Let $c$ be a section of $\C_{T^*X}$. We denote by $\delta_c$ the automorphism 
of $\C$-algebras on $T^*X$
\eq\label{eq:defdeltaa}
\delta_c:\oim{\rho}\shed[X\times \C] & \to &\oim{\rho}\shed[X\times \C] \\
P & \mapsto & \exp(c\partial_t) \circ P \circ \exp(-c\partial_t).
\nonumber
\eneq
\end{definition}

In a local coordinates system $(x,t;\xi,\tau)$, 
$$\delta_c (P(x,t,\partial_x,\partial_t)) = P(x,t+c,\partial_x,\partial_t).$$
Note that $\delta_c$ induces the identity on $\shwd[X]$.

\begin{lemma}\label{lemma:tauQCT3}
Let $\Phi:\oim{\rho}\shed[X\times \C]\to\oim{\rho}\shed[X\times \C]$
be a QCT above the identity on $U\subset T^*X$ 
satisfying the commutation properties of
Lemma \ref{lemma:tauQCT2}. Then there exist $c\in\C$ and
a section  
$P\in \shwd[X]$ of order $0$ with $\sigma_0(P)=1$
and satisfying $PP^*=1$ such that 
\eq\label{eq:qctdelta}
&&\Phi=\delta_c\circ\Ad(P).
\eneq
Moreover, such a pair $(c,P)$ is unique.
\end{lemma}
\begin{proof}
(i) In a local coordinate system, the $\tau$-preserving contact 
transformation $\psi$ above the identity on $T^*X$ is given by
\eqn
&&\psi(x,t;\xi,\tau)=(x,t+c;\xi,\tau).
\eneqn
for a locally constant function $c$.
The transformation $\Phi^{\dagger}=\Phi\circ \delta_{-c}$ is thus a QCT above 
the identity on $\opb\rho (U)\subset\dot{P}^*(X\times\C)$ commuting with $*$ 
and preserving $\partial_t$. Therefore, by Lemma \ref{le:auts}, there exists a 
unique section $P$ of $\oim{\rho}\shed[X\times \C]$ of order $0$  satisfying 
$PP^*=1$ such that $\Phi^{\dagger}=\Ad(P)$. Since $\Ad(P)(\partial_t) = 
\partial_t$, $P$ does not depend on $t$, {\em i.e.} $P$ is a section of 
$\shwd[X]$.

\noindent
(ii) Assume there exist $P$ as above, not depending on $t$ and $c\in\C$ such
that $\Ad(P)=\delta_c$. Then 
$P=\exp(c\tau)$. Hence $P$ cannot belong to $\shed[X\times \C]$, except if 
$c=0$.
\end{proof}

Denote by $\sha ut_{*, \partial_t} (\oim{\rho}\shed[X\times \C])$ the sheaf of 
QCT's above the identity on $T^*X$ commuting with $*$ and preserving 
$\partial_t$ and by $\sha ut_* (\shwd[X])$ the sheaf of QST's
 above the identity 
on $T^*X$ commuting with $*$. Denote by $\delta$ the morphism of groups 
$\C_{T^*X}\to\sha ut_{*,\partial_t}(\oim{\rho}\shed[X\times \C])$ given by 
$c\mapsto \delta_c$ and set 
\eq
&&(\shwd[X])_* = \{P\in\shwd[X];P\text{ has order }0, \sigma_0(P)=1, PP^*=1\}
,\label{eq:starcenter1}\\
&&\cor_* = \{a\in \cor ; a\text{ has order }0,\sigma_0(a)=1, 
a(\tau)a(-\tau)=1\} =(\shwd[\rmpt])_* .\label{eq:starcenter2}
\eneq
Then the diagram below commutes and has exact rows
\eqn
\xymatrix{ & 1 \ar[r] & {\C_{T^*X}\times (\shwd[X])_*} \ar[d] 
\ar[r]^-{\delta\circ\Ad} & 
{\sha ut_{*,\partial_t}(\oim{\rho}\shed[X\times \C])} \ar[d] \ar[r] & 1 \\
1 \ar[r] & {(\cor_*)_{T^*X}} \ar[r] & {(\shwd[X])_*} 
\ar[r]^-{\Ad} & {\sha ut_* (\shwd[X])} \ar[r] & 1 .}
\eneqn

\begin{remark}\label{rem:ehat3}
The results of this section still hold when replacing $\shw$ with $\shhw$.
\end{remark}

\section{Quantization of complex symplectic manifolds}\label{section:QCSM}

\begin{theorem}\label{th:qcsm}
Let $\stx$ be a complex symplectic manifold.
There exists canonically a $\cor$-abelian stack 
$\mds[\shdw,\stx]$ on $\stx$ such that if $U\subset\stx$ is an open subset 
isomorphic by a symplectic transformation $\phi$ to an open subset 
$U_X\subset T^*X$, then $\mds[\shdw,\stx]\vert_U$ is equivalent by $\phi$ to 
the stack $\mds[{\shwd[X]\vert_{U_X}}]$.
\end{theorem}

\begin{definition}
We call $\mds[\shdw,\stx]$ the stack of WKB-differential modules on $\stx$.
\end{definition}

\begin{proof}
There exists an open covering $\shu=\{U_i\}_{i\in I}$ of $\stx$ and, for each
$U_i\in\shu$, a symplectic open embedding $\phi_i:U_i\hookrightarrow T^*X_i$ 
for some cotangent bundle $T^*X_i$. 
Consider the symplectic isomorphism $\phi_{ij}=\phi_i\circ \opb \phi_j : 
\phi_j(U_{ij})\isoto \phi_i(U_{ij})$. After shrinking the covering $\shu$,
by Lemma \ref{lemma:geomconttrans} there exists a $\tau$-preserving contact 
isomorphism $\psi_{ij}$ making the diagram below commutative:
\eqn
\xymatrix{
{T^*X_j\supset \phi_j(U_{ij})}\ar[r]^-{\phi_{ij}} & {\phi_i(U_{ij}) 
 \subset T^*X_i}\\
{\dot{P}^*(X_j\times \C)\supset \opb{\rho_j}(\phi_j(U_{ij}))} \ar[u]_-{\rho_j}
 \ar[r]^-{\psi_{ij}} & {\opb{\rho_i}(\phi_i(U_{ij}))\subset
  \dot{P}^*(X_i\times\C)} \ar[u]_-{\rho_i} . }
\eneqn

\noindent
We use the notations of Lemma \ref{lemma:tauQCT2}. 
Set $\oim\rho\shed[U_i]=\opb{\phi_i}\oim{\rho_i}\shed[X_i\times\C]$ and 
$\shwd[U_i]=\opb{\phi_i}\shwd[X_i]$. Hence, $\shwd[U_i]$
is the subsheaf of $\C$-algebras of $\oim\rho\shed[U_i]$ consisting of sections 
commuting with $\partial_t$.
By shrinking the covering $\shu$ again, we may assume by Lemma 
\ref{lemma:tauQCT2} that there exist QCT's above the $\psi_{ij}$'s and hence 
isomorphisms of $\C$-algebras
\eqn
&& \Phi_{ij}:\oim\rho\shed[U_j]\vert_{U_{ij}}\isoto 
\oim\rho\shed[U_i]\vert_{U_{ij}},
\eneqn
these isomorphisms commuting with $\partial_t$ and $*$.

\noindent
With the notations of Section \ref{section:stacks} for $X=\stx$,
set $\stx_0=\bigsqcup_i U_i$, $\stx_1=\stx_0\times_{\stx}\stx_0$, etc. and let 
$j_{U_i}:U_i\hookrightarrow \stx_0$ be the natural map. Then
$\oim\rho\shed[]=\oplus_{i\in I}\eim{{j_{U_i}}}\oim\rho\shed[U_i]$ is a sheaf 
of central $\C$-algebras on $\stx_0$ endowed with an anti-involution and
$\shdw=\oplus_{i\in I}\eim{{j_{U_i}}}\shwd[U_i]$ is a sub-sheaf of algebras.
The $\Phi_{ij}$'s induce a $\C$-algebra isomorphism
$\Phi:\oim\rho\shed[1]\isoto\oim\rho\shed[0]$ commuting with $\partial_t$
and with the anti-involution.
 

\noindent
By Lemma \ref{lemma:tauQCT3}, there exist $P\in\sect(\stx_2;\shwd)$ of 
order $0$ with $\sigma_0(P)=1$, $PP^*=1$ and $c\in\sect(\stx_2;\C_{\stx_2})$ 
such that on $\stx_2$
\eqn
\Phi_{01}\circ \Phi_{12}=\Ad(P)\circ\delta_{c}\circ \Phi_{02}.
\eneqn
Now we make a computation similar to that leading to Definition 
\ref{def:lien1}.

\noindent
On $\stx_3$ one has 
\eqn
&&({\Phi}_{01}\circ {\Phi}_{12})\circ {\Phi}_{23}
   =\Ad(P_{012})\circ\delta_{c_{012}}\circ {\Phi}_{02}\circ {\Phi}_{23}\\
&&=\Ad(P_{012}P_{023})\circ\delta_{c_{012}}\circ\delta_{c_{023}}
        \circ{\Phi}_{03},
\eneqn
and 
\eqn
{\Phi}_{01}\circ ({\Phi}_{12}\circ{\Phi}_{23}) &=&{\Phi}_{01}
\circ \Ad(P_{123})\circ\delta_{c_{123}}\circ {\Phi}_{13} \\
&= &\Ad({\Phi}_{01}(P_{123}))\circ\delta_{c_{123}}
\circ {\Phi}_{01}\circ {\Phi}_{13} \\
&=& 
\Ad({\Phi}_{01}(P_{123})P_{013})\circ\delta_{c_{123}}\circ\delta_{c_{013}}
\circ {\Phi}_{03}.
\eneqn
Here, we have used the fact that the isomorphisms $\Phi$ commutes with 
$\delta_c$. It follows that 
\eq\label{eq:cocycleaijk}
&&\left\{ \begin{array}{l}
P_{012}P_{023} = {\Phi}_{01}(P_{123})P_{013} \quad \mbox{ in }
\sect(\stx_3;{\shwd}^\times),\\
c_{012}+c_{023}=c_{123}+c_{013} \quad \mbox{ in }
\sect(\stx_3;\C_{\stx_3}).
\end{array}\right.
\eneq
The $\C$-algebra isomorphism $\Phi:\oim\rho\shed[1]\isoto\oim\rho\shed[0]$ 
induces an isomorphism of $\cor$-algebras (which we denote by the same symbol):
\eqn
&& \Phi:\shwd[1]\isoto \shwd.
\eneqn
Since the isomorphisms $\delta_c$ induce the identity on $\shwd$, we get
\eqn
&&\left\{ \begin{array}{l}
\Phi_{01}\circ \Phi_{12}=\Ad(P)\circ \Phi_{02} \mbox{ on } \stx_2\\
P_{012}P_{023} = {\Phi}_{01}(P_{123})P_{013} \quad \mbox{ in }
\sect(\stx_3;{\shwd}^\times).\\
\end{array}\right.
\eneqn
Hence $\lien{W}:=(\shdw,\Phi, P)$ is an effective $\cor$-lien on 
$\stx_{\bullet}\to\stx$, and it remains to apply Theorem \ref{th:elienK}.

To prove the canonicity of the stack $\mds[\shdw,\stx]$, one argues as in the proof of Theorem \ref{th:qccm}. More precisely, for each $U_i\in \shu$ consider
another sympectic open embedding $\phi'_i:U_i\hookrightarrow T^*X'_i$. 
Setting $\oim\rho{\shed[U_i]}'=\opb{{\phi'}_i}\oim{\rho_i}\shed[X'_i\times\C]$
and ${\shwd[U_i]}'=\opb{{\phi'}_i}\shwd[X'_i]$ and proceeding as above,
we get a $\C$-algebra isomorphism 
$\Phi':\oim\rho{\shed[1]}'\isoto\oim\rho{\shed[0]}'$
commuting with $\partial_t$ and with the anti-involution and thus an
effective $\cor$-lien $\lien{W}':=({\shdw}',\Phi', P')$ on 
$\stx_{\bullet}\to\stx$ and a 2-cocycle $c'\in\sect(\stx_2;\C_{\stx_2})$. 

\noindent
By Lemma \ref{lemma:geomconttrans}, there exists a $\tau$-preserving contact 
isomorphism $\upsilon_i$ making the corresponding diagram $\eqref{diag:2}$ 
commutative. 
By Lemma \ref{le:adpdc}, there exist QCT's above the $\upsilon_i$'s and hence 
an isomorphism of $\C$-algebras
\eq
&&\Upsilon:\oim \rho\shed[]\isoto\oim\rho{\shde}'
\eneq
commuting with $\partial_t$ and $*$.
By Lemma \ref{lemma:tauQCT3}, there exist a unique pair $(Q,b)$, where
$Q\in\sect(\stx_1;\shwd)$ of order $0$ with $\sigma_0(Q)=1$, $QQ^*=1$ and 
$b\in\sect(\stx_1;\C_{\stx_1})$ are such that
\eqn
\Phi'\circ \Upsilon_1=\Ad(Q)\circ\delta_{b}\circ \Upsilon_0\circ \Phi .
\eneqn
After some computations and using the unicity of the pair $(Q,b)$, we get
\eqn
&&\left\{ \begin{array}{l}
 P'Q_{02} =  \Phi'_{01}(Q_{12})Q_{01} \Upsilon_0(P) \quad \mbox{ in }
\sect(\stx_2;{{\shwd}'}^\times),\\
c'+b_{02}=b_{12}+ b_{01} +c \quad \mbox{ in }
\sect(\stx_2;\C_{\stx_2}).
\end{array}\right.
\eneqn
Hence the 2-cocycle $c$ and $c'$ are cohomologous and the pair $(\Phi,Q)$ 
defines an effective isomorphism of $\cor$-liens $\lien{W}\isoto\lien{W}'$. 
\end{proof}

\begin{remarks}
(i) The Cech cohomology class in $H^2(\stx;\C)$ defined by the 2-cocycle 
$c$ in \eqref{eq:cocycleaijk} coincides with the class $-[\omega]$ 
of the symplectic form $\omega$ on $\stx$ under the de Rham
isomorphism. 
Indeed, let $\shu =\{U_i\}_{i\in I}$ be the open covering of $\stx$ as in the 
proof of Theorm \ref{th:qcsm} and let $\alpha_{T^*X_i}$ be the canonical 
1-form on $T^*X_i$. Then $\omega \vert_{U_i}=d (\phi^*_i\alpha_{T^*X_i})$ 
and $(\phi^*_j\alpha_{T^*X_j})\vert_{U_{ij}} = (\phi^*_i\alpha_{T^*X_i})
\vert_{U_{ij}} + da_{ij}$, for a function $a_{ij}$ on $U_{ij}$.
By Lemma \ref{lemma:geomconttrans}, the contact isomorphisms $\psi_{ij}$ in 
\eqref{diag:2} are given by
$$
(p,(t;\tau))\mapsto (\phi_{ij}(p \opb\tau),(t+ a_{ij}(p \opb\tau);\tau)).
$$ 
Hence the $c_{ijk}$ coincide with $a_{ij}+a_{jk}-a_{ik}$, {\em i.e.} the Cech
class given by the ${c_{ijk}}$'s corresponds to the class $-[\omega]$ under 
the de Rham isomorphism.

\noindent
(ii) 
Recall that the group $\cor _*$ has been defined in \eqref{eq:starcenter2}.
For each 2-cocyle $ a\in \sect(\stx_2;(\cor _*)_{\stx})$ we may associate 
a $\cor$-abelian stack $\mds[\shdw,a,\stx]$ on $\stx$ which is locally 
equivalent to the stack $\mds[\shdw,\stx]$. Indeed,
$\mds[\shdw,a,\stx]$ is the stack associated with the effective lien 
$\lien{W}_a:=(\shdw,\Phi, aP)$.  It is straightforward to check 
that, if $b$ is another 2-cocyle cohomologous to $a$, there is an 
effective isomorphism of $\cor$-liens $\lien{W}_a\isoto\lien{W}_b$ and then
an equivalence of  $\cor$-stacks $\mds[\shdw,a,\stx]\simeq \mds[\shdw,b,\stx]$.

\noindent
We refer to \cite{De} for a similar construction in the framework of
deformation algebras on real symplectic manifolds.
\end{remarks}

\begin{remark}\label{rem:kons}
One can extend the results of Section \ref{section:Invquant} 
to the symplectic case and treat
smooth involutive submanifolds of complex symplectic manifolds. 
But one should recall that
Kontsevich \cite{Ko} has obtained a much more general result, quantizing
arbitrary complex Poisson manifolds.
\end{remark}
\begin{remark}\label{rem:ehat4}
The results of this section still hold when replacing $\shw$ with $\shhw$.
\end{remark}

\vspace*{1cm}

\parbox[t]{15em}
{\scriptsize{
\noindent
Pietro Polesello\\
Universit{\`a} di Padova\\
Dipartimento di Matematica\\
via G. Belzoni, 7
35131 Padova, Italy\\
or: Universit{\'e} Pierre et Marie Curie\\
Institut de Math{\'e}matiques\\
175, rue du Chevaleret, 75013 Paris France\\
pietro@math.jussieu.fr\\}}
\qquad
\parbox[t]{15em}
{\scriptsize{
Pierre Schapira\\
Universit{\'e} Pierre et Marie Curie\\
Institut de Math{\'e}matiques\\
175, rue du Chevaleret, 75013 Paris France\\
schapira@math.jussieu.fr\\
http://www.math.jussieu.fr/{\~{}} schapira/}}

\end{document}